\documentclass[12pt]{amsart}
\usepackage{amssymb, hyperref}

\newcommand{\Bgp}{{\Z^\N}}

\long\def\forget#1\forgotten{}
\newcommand{\issuenumber}{33}
\newcommand{\issuemonth}{February}
\newcommand{\issueyear}{2012}

\newcommand{\ed}{
\newpage

\section{Unsolved problems from earlier issues}

\begin{issue}Is $\binom{\Omega}{\Gamma}=\binom{\Omega}{\Tau}$?\end{issue}%\stepcounter{issue}
\begin{issue}Is $\ufin(\cO,\Omega)=\sfin(\Gamma,\Omega)$?And if not, does $\ufin(\cO,\Gamma)$ imply
$\sfin(\Gamma,\Omega)$?\end{issue}%\stepcounter{issue}
\stepcounter{issue}\begin{issue}Does $\sone(\Omega,\Tau)$ imply $\ufin(\Gamma,\Gamma)$?\end{issue}
\begin{issue}Is $\fp=\fp^*$? (See the definition of $\fp^*$ in that issue.)\end{issue}
\begin{issue}Does there exist (in ZFC) an uncountable set satisfying $\sfin(\cB,\cB)$?\end{issue}
\stepcounter{issue}
\begin{issue}Does $X \nin \NON(\cM)$ and $Y\nin\mathsf{D}$ imply that $X\cup Y\nin \COF(\cM)$?\end{issue}
\begin{issue}[CH]Is $\split(\Lambda,\Lambda)$ preserved under finite unions?\end{issue}
\begin{issue}Is $\cov(\cM)=\fo$? (See the definition of $\fo$ in that issue.)\end{issue}
\stepcounter{issue}
\begin{issue}Could there be a Baire metric space $M$ of weight $\aleph_1$ and a partition
$\mathcal{U}$ of $M$ into $\aleph_1$ meager sets where for each ${\mathcal U}'\subset\mathcal U$,
$\bigcup {\mathcal U}'$ has the Baire property in $M$?\end{issue}
\stepcounter{issue} %% no problem in Issue 13
\begin{issue}Does there exist (in ZFC) a set of reals $X$ of cardinality $\fd$ such that all
finite powers of $X$ have Menger's property $\sfin(\cO,\cO)$?\end{issue}
\begin{issue}Can a Borel non-$\sigma$-compact group be generated by a Hurewicz subspace?\end{issue}
\begin{issue}[MA]Is there $X\sbst\bbR$ of cardinality continuum, satisfying $\sone(\BO,\BG)$?\end{issue}
\begin{issue}[CH]Is there a totally imperfect $X$ satisfying $\ufin(\cO,\Gamma)$
that can be mapped continuously onto $\Cantor$?\end{issue}
\begin{issue}[CH]Is there a Hurewicz $X$ such that $X^2$ is Menger but not Hurewicz?\end{issue}
\begin{issue}Does the Pytkeev property of $C_p(X)$ imply that $X$ has Menger's property?\end{issue}
\begin{issue}Does every hereditarily Hurewicz space satisfy $\sone(\BG,\BG)$?\end{issue}
\begin{issue}[CH]Is there a Rothberger-bounded $G\le\Bgp$ such that $G^2$ is not Menger-bounded?\end{issue}
\begin{issue}Let $\cW$ be the van der Waerden ideal. Are $\cW$-ultrafilters closed under products?\end{issue}
\begin{issue}Is the $\delta$-property equivalent to the $\gamma$-property $\binom{\Omega}{\Gamma}$?\end{issue}
\stepcounter{issue}\stepcounter{issue}
\general\end{document}}

\newcommand{\Cantor}{{\{0,1\}^\N}}

\newcommand{\fd}{\mathfrak{d}}
\newcommand{\fp}{\mathfrak{p}}
\newcommand{\NON}{{\mathsf   {NON}}}\newcommand{\COF}{{\mathsf   {COF}}}

\newcommand{\cM}{\mathcal{M}}
\newcommand{\cov}{\mathsf{cov}}

\newcommand{\bbR}{\mathbb{R}}
\newcommand{\fo}{\mathfrak{od}}

\renewcommand{\split}{\mathsf{Split}}\newcommand{\bq}{\begin{quote}}\newcommand{\eq}{\end{quote}}
\newcommand{\cO}{\mathcal{O}}\newcommand{\cB}{\mathcal{B}}\newcommand{\BG}{\cB_\Gamma}
\newcommand{\BO}{\cB_\Omega}

\newcommand{\sone}{\mathsf{S}_1}\newcommand{\sfin}{\mathsf{S}_\mathrm{fin}}
\newcommand{\ufin}{\mathsf{U}_\mathrm{fin}} 
\newcommand{\nin}{\not\in}
\newcommand{\cW}{\mathcal{W}}

\newcommand{\N}{\mathbb{N}}\newcommand{\Z}{\mathbb{Z}}
\newcommand{\sbst}{\subseteq}
\newcommand{\by}[2]{\par\hfill\emph{#1}, #2}\newcommand{\nby}[1]{\par\hfill\emph{#1}}\newcommand{\Tau}{\mathrm{T}}
\newcommand{\CE}{\textsc{CE}}
\newtheorem{thm}{Theorem}[section]\newcommand{\bthm}{\begin{thm}} \newcommand{\ethm}{\end{thm}}
\newtheorem{prop}[thm]{Proposition}\newcommand{\bprp}{\begin{prop}} \newcommand{\eprp}{\end{prop}}
\newtheorem{fact}[thm]{Fact}\newcommand{\bfct}{\begin{fact}} \newcommand{\efct}{\end{fact}}
\newtheorem{prob}[thm]{Problem}\newcommand{\bprb}{\begin{prob}} \newcommand{\eprb}{\end{prob}}
\newtheorem{lem}[thm]{Lemma}\newcommand{\blem}{\begin{lem}} \newcommand{\elem}{\end{lem}}
\newtheorem{claim}[thm]{Claim}\newcommand{\bclm}{\begin{claim}} \newcommand{\eclm}{\end{claim}}
\newtheorem{cor}[thm]{Corollary}\newcommand{\bcor}{\begin{cor}} \newcommand{\ecor}{\end{cor}}
\newtheorem{conj}[thm]{Conjecture}\newcommand{\bcnj}{\begin{conj}} \newcommand{\ecnj}{\end{conj}}
\theoremstyle{definition}\newtheorem{defn}[thm]{Definition}\newcommand{\bdfn}{\begin{defn}} \newcommand{\edfn}{\end{defn}}
\theoremstyle{remark}\newtheorem{rem}[thm]{Remark}\newcommand{\brem}{\begin{rem}} \newcommand{\erem}{\end{rem}}
\newtheorem{cnv}[thm]{Convention}\newcommand{\bcnv}{\begin{cnv}} \newcommand{\ecnv}{\end{cnv}}
\newtheorem{exam}[thm]{Example}\newcommand{\bexm}{\begin{exam}} \newcommand{\eexm}{\end{exam}}
\newtheorem{issue}{Issue}\newcommand{\bpf}{\begin{proof}} \newcommand{\epf}{\end{proof}}
\newcommand{\be}{\begin{enumerate}}\newcommand{\ee}{\end{enumerate}}\newcommand{\bi}{\begin{itemize}}
\newcommand{\ei}{\end{itemize}}
\newcommand{\general}{\small\vfill\par\noindent\hrulefill\par
\noindent\textbf{Previous issues.} The previous issues of this
bulletin are available online at\\
\url{http://front.math.ucdavis.edu/search?\&t=\%22SPM+Bulletin\%22}
\\[0.1cm]
\textbf{Contributions.} Announcements, discussions, and open problems should be emailed
to \texttt{tsaban@math.biu.ac.il}\\[0.1cm]
\textbf{Subscription.}
To receive this bulletin (free) to your e-mailbox, e-mail us.
}

\newcommand{\arXivl}[4]{\subsection{#2}{#4}\par\hfill{\arx{#1}}\par\hfill\emph{#3}}
\newcommand{\arXiv}[3]{\subsection{#2}\mbox{}\par\hfill{\arx{#1}}\par\hfill\emph{#3}}

\newcommand{\arx}[1]{\url{http://arxiv.org/abs/#1}}

\title[$\mathcal{SPM}$ Bulletin \textbf{\issuenumber} (\issuemonth{} \issueyear)]{%
$\mathcal{SPM}$ Bulletin\\[0.5cm]
Issue number \issuenumber: \issuemonth{} \issueyear{} \CE{}}

\begin{document}
\maketitle

%\tableofcontents

\section{Editor's note}

Dear Friends,

You may have heard by now of the tragic loss of Ireneusz Rec\l{}aw. 
A lecture of Dr.\ Rec\l{}aw, which I heard as an an undergraduate student at Bar-Ilan University, has set
my interest in special sets of reals. This issue is dedicated to his memory.

Jakub Jasi\'nski provides, kindly, a memorandum below.

\medskip

\by{Boaz Tsaban}{tsaban@math.biu.ac.il}

\hfill \texttt{http://www.cs.biu.ac.il/\~{}tsaban}

\subsection*{Irek Rec\l{}aw}
Ireneusz (Irek) was my dear friend for over 34 years. I had the privilege to be his teacher, a colleague and a coauthor. Over the years, frequently visiting each other across the Atlantic we and our families developed a unique form of friendship.

Dr.\ Ireneusz Reclaw was a professor of Mathematics at the University of Gdansk. He shared his love of mathematics with his students in Gdansk Poland, Scranton Pennsylvania, Auburn Alabama and Berlin Germany. Well respected in the world community of experts in set and real function theory he authored over 40 publications in the most reputable math journals. For many years Ireneusz suffered from a chronic kidney disease, diagnosed with gallbladder cancer in December 2011, died February 04, 2012 at the age of 51.

His uncommon intellect not only lead him to numerous mathematical discoveries but he will also be remembered for his warm personality, witty sense of humor, love of people, understanding of  nature and knowledge of  birds.

\nby{Jakub Jasi\'nski}

\section{Long announcements}

\arXivl{1112.2373}
{Pointwise convergence of partial functions: The Gerlits-Nagy Problem}
{Tal Orenshtein, Boaz Tsaban}
{For a set $X\sbst\bbR$, let $B(X)\sbst\bbR^X$ denote the space of Borel real-valued functions on $X$,
with the topology inherited from the Tychonoff product $\bbR^X$.
Assume that for each countable $A\sbst B(X)$, each $f$ in the closure of
$A$ is in the closure of $A$ under pointwise limits of sequences of partial functions.
We show that in this case, $B(X)$ is countably Fr\'echet-Urysohn, that is,
each point in the closure of a countable set is a limit of a sequence of elements of that set.
This solves a problem of Arnold Miller. The continuous version of this problem
is equivalent to a notorious open problem of Gerlits and Nagy.
Answering a question of Salvador Herna\'ndez, we show that
the same result holds for the space of all Baire class $1$ functions on $X$.

We conjecture that the answer to the
continuous version of this problem is negative, but we identify a nontrivial class of sets
$X\sbst\bbR$, for which we can provide a positive solution to this problem.

The proofs establish new local-to-global correspondences,
and use methods of infinite-combinatorial topology,
including a new fusion result of Francis Jordan.}

\arXivl{1201.1576}
{Selected results on selection principles}
{Ljubisa D.R. Kocinac}
{We review some selected recent results concerning selection principles in
topology and their relations with several topological constructions.

Published in: Proceedings of the Third Seminar on Geometry and Topology, July
 15--17, 2004, Tabriz, Iran, pp. 71--104.
}

\arXivl{1201.4909}
{Remarks on countable tightness}
{Marion Scheepers}
{Countable tightness may be destroyed by countably closed forcing. We
characterize the indestructibility of countable tightness under countably
closed forcing by combinatorial statements similar to the ones Tall used to
characterize indestructibility of the Lindel\"of property under countably closed
forcing. We consider the behavior of countable tightness in generic extensions
obtained by adding Cohen reals. We show that HFD's are indestructibly countably
tight.}

\arXivl{1202.0194}
{Weak covering properties and infinite games}
{Liljana Babinkostova, Bruno A. Pansera and Marion Scheepers}
{We investigate game-theoretic properties of selection principles related to
weaker forms of the Menger and Rothberger properties. We prove theorems
characterizing the existence of winning strategies for player TWO in several of
these games. We show that for appropriate spaces several of these selection
principles are characterized in terms of a corresponding game. We use generic
extensions by Cohen reals and certain preservation results to discuss the
necessity of some of the hypotheses in our theorems.}

\arXivl{1202.1619}
{Some Open Problems in Topological Algebra}
{Taras Banakh, Mitrofan Choban, Igor Guran, Igor Protasov}
{This is the list of open problems in topological algebra posed on the
conference dedicated to the 20th anniversary of the Chair of Algebra and
Topology of Lviv National University, that was held on 28 September 2001.

Journal reference: Visnyk Lviv Univ.\ 61 (2003), 13--20.
}

\section{Short announcements}\label{RA}

\arXiv{1109.6517}
{Local properties on the remainders of the topological groups}
{Fucai Lin}

\arXiv{1110.4153}
{Tukey types of ultrafilters}
{Natasha Dobrinen and Stevo Todorcevic}

\arXiv{1110.4154}{Continuous cofinal maps on ultrafilters}
{Natasha Dobrinen}

\arXiv{1110.4483}
{Approximable WAP- and LUC-interpolation sets}
{Jorge Galindo and Mahmoud Filali}

\arXiv{1111.6705}
{A Ramsey-Classification Theorem and its Application in the Tukey Theory
 of Ultrafilters}
{Authors: Natasha Dobrinen, Stevo Todorcevic}

\arXiv{1201.1568}
{The C-compact-open topology on function spaces}
{Alexander V. Osipov}

\arXiv{1201.3814}
{The weights of closed subgroups of a locally compact group}
{Salvador Hern\'andez, Karl H. Hofmann, and Sidney A. Morris}

\arXiv{1201.5181}
{The well-ordered (F) spaces are D-spaces}
{Xu Yuming}

\ed